\documentclass[12pt]{article}
\usepackage{amsfonts}

\usepackage{amsmath}
\usepackage{epsfig}
\usepackage{amsthm}
\usepackage{color}

\font\tenopen=msbm10 \font\sevenopen=msbm7 \font\fiveopen=msbm5
\newfam\openfam
\def\open{\fam\openfam\tenopen}
\textfont\openfam=\tenopen
\scriptfont\openfam=\sevenopen
\scriptscriptfont\openfam=\fiveopen

\def\R{{\open \mathbb{R}}}

\def\E {{\open \mathbb{E}}}

\newtheorem{theorem}{Theorem}

\newtheorem{corollary}[theorem]{Corollary}

\newtheorem{lemma}[theorem]{Lemma}

\newtheorem{property}[theorem]{Property}
\newtheorem{proposition}[theorem]{Proposition}

\newcommand{\bb}{\begin{eqnarray*}}

\newcommand{\ee}{\end{eqnarray*}}

\input{tcilatex}

\begin{document}

\title{}

\begin{center}
{\Large \textbf{Another approach to Brownian motion}}\\[0pt]
Magda Peligrad${}^{a}$\footnote{\textit{Mathematical Subject Classification}
(2000):60G51, 60F05
\par
Supported in part by a Charles Phelps Taft Memorial Fund grant and a NSA\
grant.} \textit{and } Sergey Utev${}^{b}$
\end{center}

\vskip5pt \noindent ${}^{a}$Department of Mathematical Sciences, University
of Cincinnati, PO Box 210025, Cincinnati, OH 45221-0025, USA \vskip5pt

\noindent ${}^{b}$School of Mathematical Sciences, University of Nottingham,
Nottingham, NG7 2RD, UK

\vskip10pt \noindent \textbf{Abstract}\vskip10pt

Motivated by the central limit theorem for weakly dependent variables, we
show that the Brownian motion $\{X(t);t\in \lbrack 0,1]\}$, can be modeled
as a process with independent increments, satisfying the following limiting
condition 
\begin{equation*}
\lim \inf_{h\downarrow 0}\E f(h^{-1/2}[X(s+h)-X(s)])\;\geq \;\E f(X(1))
\end{equation*}
almost surely for all$\;0\leq s<1$, where $\E f(X(1))<\infty $ and $f:%
\mathbb{R}\rightarrow \mathbb{R}$ is a symmetric, continuous, convex
function with $f(0)=0$, strictly increasing on $\mathbb{R}^{+}$ and
satisfying the following growth condition: 
\begin{equation*}
f(Kx)\leq K^{p}f(x),\text{for a certain }p\in \lbrack 1,2),\text{all }K\geq
K_{0}\text{ and all }x>0\;
\end{equation*}
(for example, $f(x)=x^{p}[A+B\ln (1+Cx)]$, with $x>0,$ $p\in \lbrack 1,2)$, $%
A>0$ and $B,C\geq 0$).\newline

\bigskip

\bigskip

\textit{Key words}: Levy process, Brownian motion, processes with
independent increments, central limit theorem, weakly dependent sequences.

\section{Introduction}

\quad \quad For the partial sums $S_{n}=Y_{1}+\ldots +Y_{n}$ of a centered
stationary strongly mixing sequence $\{Y_{i}\}$ with finite second moment,
the well-known sufficient conditions for the central limit theorem are that $%
\mathrm{Var}(S_{n})/n$ is slowly varying as $n\rightarrow \infty $ and the
sequence $\{S_{n}^{2}/\sigma _{n}^{2}\}$ is uniformly integrable $(\sigma
_{n}^{2}=\mathrm{Var}(S_{n})).$ The conditions are checkable under various
mixing conditions and they lead to the central limit theorem under the
normalization $\sigma _{n}$(see, Denker (1986) and Peligrad (1986) for a
survey).

Dehling, Denker and Philipp (1986) proved an interesting central limit
theorem using the non-traditional normalization $\rho _{n}=\E|S_{n}|$. One
of their results, Theorem 3, roughly states that if both sequences $\sigma
_{n}^{2}/n$ and $\rho _{n}/\sqrt{n}$ are slowly varying as $n\rightarrow
\infty $, then the central limit theorem holds.

On the other hand, Braverman, Mallows and Shepp (1995) showed that if the
absolute moments of partial sums of i.i.d. symmetric variables are equal to
those of normal variables, then the marginals have normal distribution. This
fact suggested the conjecture that probably the absolute moments alone
characterize the homogeneous process with independent increments (see Bryc
(2002) for a discussion on this topic and related conjectures).

Our main interest in this topic is to prove some of these conjectures and to
apply them to understand the nature of the intricate normalization in
Dehling, Denker and Philipp (1986).

Throughout the paper, $\{W(t);t\in \lbrack 0,1]\}$ denotes the standard
Brownian motion, i.e. a Gaussian process $\{W(t);t\in \lbrack 0,1]\}$ with
independent increments, $\E[W(t)]=0$ and $\E[W(t)W(s)]=\min (t,s)$. By $W$
we denote a standard normal variable. Also $\mu $ denotes the Lebesgue
measure, $h\downarrow 0$ denotes convergence over positive real numbers, $%
[x] $ denotes the integer part of $x$. For two processes with independent
increments $\{X(t)\;;\;[0,1]\}$ and $\{Y(t);t\in \lbrack 0,1]\}$, equality $%
X(t)=Y(t)$ means that their increments have the same distribution.

The process $\{X(t);t\in \lbrack 0,1]\}$, is called homogeneous if $%
X(t+s)-X(t)=^{d}X(s)$ where $=^{d}$ means equality in distribution. Finally,
the process $\{X(t);t\in \lbrack 0,1]\}$ is called stochastically continuous
if it does not have deterministic jumps, i.e. $P(|X(t+s)-X(t)|>u)\rightarrow
0$ as $s\rightarrow 0$ for any $u>0$ and $t\in \lbrack 0,1]$.

Our paper is organized as follows. In Section 2 we include the
representation results and their corollaries. Section 3 is dedicated to
their proofs. In section 4 we give an application of the characterization
results to the central limit theorem.

\section{Characterization Results}

As a class of potential characterizing functions, we consider non-negative
functions satisfying the following conditions: 
\begin{eqnarray}
&&%
\mbox{The function $f:\R\to\R$ is symmetric, continuous, convex,
strictly increasing on $\R^+$, $f(0)=0$,}  \notag \\
&&%
\mbox{and there exists $p\in[1,2)$ and $K_0\geq 0$ such that $f(Kx)\leq K^p f(x)$
for all $K\geq K_0$}\;.  \label{fcond}
\end{eqnarray}
For example, $f(x)=x^{p}$, or more generally, $f(x)=x^{p}[A+B\ln (1+Cx)]$, $%
x>0$, for a $p\in \lbrack 1,2),$ $A>0$ and $B,C\geq 0$ satisfies (\ref{fcond}%
) for some $p\prime $, $p<p\prime <2$ .

The following theorem is the main result of the paper.

\begin{theorem}
\label{tm1}Let $f$ be a positive function satisfying (\ref{fcond}) and let $%
\{X(t);$ $t\in \lbrack 0,1]\}$ be a process with independent increments, $%
X(0)=0,$ and $\E f(X(1))<\infty $. Assume in addition that $\mu $ - almost
surely for $s\in \lbrack 0,1]$, 
\begin{equation}
\lim \inf_{h\downarrow 0}\E f(h^{-1/2}[X(s+h)-X(s)])\;\geq \;\E f(X(1))\;.
\label{keycond}
\end{equation}
Then, $\{X(t);$ $t\in \lbrack 0,1]\}$ is a Gaussian process that admits the
representation 
\begin{equation}
X(t)=\sigma W(t)+\E X(t)\;\;\qquad \mbox{for all }\;t\in \lbrack 0,1]
\label{repres}
\end{equation}
where $\E X(1)=0$, $\sigma =\Psi ^{-1}(\E f(X(1))$ and the function $\Psi
(x)=\E f(x W)$ is continuous and strictly increasing for $x>0$.
\end{theorem}

\begin{corollary}
\label{tm1cor1}Let $f$ satisfies (\ref{fcond}) and let $\{X(t);$ $t\in
\lbrack 0,1]\}$ be a stochastic process with independent increments, $%
X(0)=0, $ satisfying the following condition: 
\begin{equation*}
\E f(t^{-1/2}[X(t+s)-X(s)])=\E f(W)\qquad \mbox{for all }\;0\leq s\leq
s+t\leq 1 \;.
\end{equation*}
Then, $\{X(t);t\in \lbrack 0,1]\}$ is a standard Brownian motion.
\end{corollary}

By taking $f(x)=x$, the corollary gives an affirmative answer to a
conjecture of Bryc and Peligrad formulated in a survey paper by Bryc (2002).

We notice that we do not impose any conditions on the sample path properties
of the stochastic process $\{X(t);$ $t\in \lbrack 0,1]\}$. In particular, a
Gaussian process satisfying (\ref{repres}) does not have to be a
semimartingale (see for example, Jacod and Shiryaev, p.106).

For a stochastically continuous homogeneous processes, it is enough to check
the limiting condition in (\ref{keycond}) only on one subsequence, which is
useful in applications.

\begin{corollary}
\label{tm1cor2}Suppose that $\{X(t);$ $t\geq 0\}$ is a stochastically
continuous homogeneous process, with independent increments, $X(0)=0$ (i.e.
Levy process), $\E f(X(1))<\infty $, and assume there exists a positive
sequence $t_{n}\rightarrow 0$ such that 
\begin{equation}
\lim \inf_{t_{n}\rightarrow 0}\E f(t_{n}^{-1/2}[X(t_{n})])\;\geq \;\E %
f(X(1))\;.  \label{keycondh}
\end{equation}
Then, $X(t)=\sigma W(t)$ for all $t\in \lbrack 0,1]$ where $\sigma $ is
defined as in Theorem \ref{tm1}.
\end{corollary}

\quad In the following proposition, we show that without the stochastic
continuity assumption in Corollary \ref{tm1cor2}, the result is not true in
general.

\begin{proposition}
\label{example}There exists a non-Gaussian homogeneous stochastic process $%
\{X(t);t\geq 0\}$ with independent increments, with $X(0)=0$, such that (\ref
{keycondh}) is satisfied with some positive sequence $t_{n}\rightarrow 0$.
\end{proposition}

We notice that the restriction $p<2$ in (\ref{fcond}) in Theorem \ref{tm1}
is necessary, in general. For example, if $f(x)=x^{p}$ with $p\geq 2$ or
more generally $f(x)$ is a bounded twice continuously differentiable
function with $f(0)=f^{\prime }(0)=0$, then (\ref{keycond}) is a condition
only on the variance of the increments $X(t+u)-X(t)$ and does not imply (\ref
{repres}).

\section{Proofs}

\quad The proof is divided in a few separate lemmas, some of them are of the
independent interest.

%\newline
In the first lemma, we state some properties of the function $f(x)$
satisfying Condition (\ref{fcond}).

\begin{lemma}
\label{fprop} Suppose that the function $f$ satisfies Condition (\ref{fcond}%
). Then,\newline
(a) There exists a positive $\alpha >0$ such that for all $t\geq 2$ and $x>0$%
, $f(tx)\leq t^{\alpha }f(x)$.\newline
(b) In addition, there exists a positive constant $C$ such that for all $%
x,y\geq 0$ and $z>1$, 
\begin{equation*}
f(x+y)\leq C[f(x)+f(y)]\;,\text{ \ }f(x)\leq C(x+x^{2})\quad \text{ }%
\mbox{and}\quad f(z)\geq z/C\;.
\end{equation*}
\end{lemma}

\textit{Proof}. To prove the statement (a), we assume without loss of
generality that $K_{0}>2$ in Condition (\ref{fcond}). Then, for $t>K_{0}$,
we know that $f(tx)\leq t^{p}f(x)$. For $2\leq t\leq K_{0}$, 
\begin{equation*}
f(tx)=f(K_{0}(tx/K_{0}))\leq K_{0}^{p}f(tx/K_{0})\leq K_{0}^{p}f(x)=2^{p\log
_{2}(K_{0})}f(x)
\end{equation*}
which proves (a) with $\alpha =p\log _{2}(K_{0})$.

First inequality in part (b) is a simple consequence of (a) since 
\begin{equation*}
f(x+y)\leq \tfrac{1}{2}(f(2x)+f(2y))\leq 2^{a-1}(f(x)+f(y))
\end{equation*}
The other two assertions are simple consequences of Condition (\ref{fcond}).$%
\diamond $

\bigskip

In the next lemma, we analyze some properties of expectations associated to
the function $f(x)$ satisfying condition (\ref{fcond}).

\begin{lemma}
\label{Glemma} Let $f$ satisfies Condition (\ref{fcond}) and let $W$ be a
standard normal random variable.

(a) Let $G(y)=\E f(W+y).$ Then the function $G$ is symmetric, continuous and
strictly increasing for $y>0$. Also, the function $\Psi (x)=\E f(x|W|)$ is
continuous and strictly increasing for $x>0$.

(b) Assume that $Y$ is a random variable independent of $W$ and let $x\geq 0$%
. Then, $\E f(xW+Y)\geq \E f(xW),$ and the equality holds if and only if $%
P(Y=0)=1$.

(c) Assume that $X$ and $Y$ are independent random variables and $\E %
f(X+Y)<\infty $. Then, $\E f(X)$ and $\E f(Y)<\infty $.

(d) Assume that $X$ and $Y$ are i.i.d. random variables with $\E(X)=0$. Then
there is a constant $C_{1}$ which depends only on $p$ from Condition (\ref
{fcond}), such that $\E f(X)\leq \E f(X-Y)\leq C_{1}\E f(X)$.
\end{lemma}

\textit{Proof.} Notice first that $G$ is infinitely differentiable. In
addition, $G$ is symmetric, since the random variable $W$ is symmetric, and $%
G$ is convex, since the function $f$ is convex. Moreover, since $W$ has as
support all the real numbers, and $f$ is non-constant, the function $G$ is
strictly convex. We shall also notice that, by symmetry, $G^{\prime }(0)=0$
and the function $G^{\prime }(x)$ is strictly positive for $x>0$. The same
argument works for the function $\Psi (x)$ which proves (a). Statement (c)
follows from the Fubini theorem, since an a.s. finite convex function is
finite. Finally, statement (d) follows from the Jensen inequality,
monotonicity of the function $f$ on $\mathbb{R}^{+}$ and Property (b) in
Lemma \ref{fprop}.$\diamond $

\bigskip

The following moment inequality was established by Klass and Nowicki (1997,
Lemma 2.6). Although, their result was stated for $A\leq 1$ the adaptation
is immediate by considering blocks with partial sums satisfying (\ref{KNcond}%
) with $A\leq 1$. We also formulate this lemma for an infinite number of
pairs by passing to the limit.

\begin{lemma}
\label{KNlm} Let $\{(X_{k},I_{B_{k}});$ $k\geq 1\}$ be independent pairs of
random variables, where $I_{B}$ is an indicator variable. Assume that the
function $H:{\fam\openfam\tenopen\mathbb{R}}\rightarrow {\fam\openfam\tenopen%
\mathbb{R}}$ is symmetric, continuous, strictly increasing on ${\fam\openfam%
\tenopen\mathbb{R}}^{+}$, $H(0)=0$ and there is a $p>0$, such that $%
H(Kx)\leq K^{p}H(x)$ for all $K\geq 2,$ $x>0$. Suppose that 
\begin{equation}
\sum_{i\geq 1}P(B_{i})\leq A\;.  \label{KNcond}
\end{equation}
Then, there exist two positive constants $c_{1}$ and $c_{2}$ such that 
\begin{equation*}
c_{1}\sum_{i\geq 1}\E H(X_{i}I_{B_{i}})\;\leq \;\E H\left( \sum_{i\geq
1}X_{i}I_{B_{i}}\right) \leq \;c_{2}\sum_{i\geq 1}\E H(X_{i}I_{B_{i}})\;.
\end{equation*}
\end{lemma}

The next property is going to be used several times in the proofs (see for
example, Rogers, (1998), namely Vitali's argument in theorems 63 and 64)).

\begin{property}
\label{pr1} Assume that $F(x)$ is a non-decreasing function on $[0,1]$, then 
\begin{equation*}
\mu \{s\in \lbrack 0,1):\lim \sup_{h\rightarrow 0}\;h^{-1}|F(s+h)-F(s)|\geq
K\}\;\leq \;(F(1)-F(0))/K\;.
\end{equation*}
\end{property}

\quad The following technical lemma is useful for handling the
non-stationary case.

\begin{lemma}
\label{funct} For any function $c(t),$ $t\in \lbrack a,b]$, 
\begin{eqnarray*}
\lim \inf_{h\downarrow 0}\;|c(s+h)-c(s)|/\sqrt{h}=0\quad \mu -%
\mbox{almost
surely} \;.
\end{eqnarray*}
\end{lemma}

\textit{Proof.} Without loss of generality we can take $[a,b]=[0,1]$ and
notice first that 
\begin{eqnarray*}
&&\left\{ s\in (0,1):\lim \inf_{h\downarrow 0}\;|c(s+h)-c(s)|/\sqrt{h}%
>0\right\} =\cup _{n,k,m=1}^{\infty }\cup _{j=0}^{k-1}A_{n,k,j,m}\quad %
\mbox{where}\;A_{n,k,j,m} \\
&=&\left\{ s\in (j/k,(j+1)/k]\;:\;|c(s)|<m\;,\;|c(s+h)-c(s)|\geq \sqrt{h}%
/n\quad \mbox{for all}\;h\in (0,1/k)\right\} \;\;.
\end{eqnarray*}
We say that the set $G\subseteq \lbrack 0,1]$ and the function $c$ satisfy
Property $(G,c)$ if there exist two positive real numbers $u$ and $w$ such
that 
\begin{equation*}
(G,c):\qquad \mbox{for all}\quad s,s_{1},s_{2}\in G,\quad
|c(s)|<w\;,\;|c(s_{2})-c(s_{1})|\geq u\sqrt{|s_{2}-s_{1}|}\;.
\end{equation*}
Clearly, the set $A_{n,k,j,m}$ and the function $c$ satisfy Property $(G,c)$
with $u=1/n$ and $w=m$ whence, it is enough to show that if $G$ and $c$
satisfy Property $(G,c)$, then $\mu (G)=0$.

Let $D=c(G),$ that is the image of the set $G$. We observe that the function 
$c:G\rightarrow D$ is one to one and let $d=c^{-1}:D\rightarrow G$ be its
inverse function. Then, the set $D$ and the function $d$ satisfy the
property 
\begin{equation*}
(D,d):\qquad \mbox{for all}\quad u_{1},u_{2}\in D,\quad
|d(u_{2})-d(u_{1})|\leq (1/u^{2})|u_{2}-u_{1}|^{2},\text{ \ }D\subset
(-w,w)\;.
\end{equation*}
Let $N$ be a positive integer, $\delta =w/N$ and define the intervals $%
\Delta _{i}=(\delta i,\delta i+\delta ]$, $i=-N,1-N,\ldots ,N-1$. Then, 
\begin{equation*}
D\subseteq \cup _{i=-N}^{N-1}D\cap \Delta _{i}\quad \mbox{and so}\quad
d(D)\subseteq \cup _{i=-N}^{N-1}d(D\cap \Delta _{i})
\end{equation*}
which implies that the outer measure $\mu ^{\ast }$ of the set $G$ is
bounded by 
\begin{equation}
\mu ^{\ast }(G)=\mu ^{\ast }(d(D))\leq \sum_{i=-N}^{N-1}\mu ^{\ast }[d(D\cap
\Delta _{i})]\;.  \label{bound1}
\end{equation}
Further, we use the following ideas associated with the computation of the
Hausdorff measurers. \newline
The first idea is a standard upper bound on the outer measure of a set by
its diameter 
\begin{equation*}
\mu ^{\ast }(A)\leq \mathrm{diam}(A):=\sup \{|x-y|\;:\;x,y\in A\}
\end{equation*}
The second idea is the bound on the diameter of the image of the Lipschitz
function $g:T\rightarrow W$, 
\begin{equation*}
\mathrm{diam}(g(T))\leq \mathrm{diam}(T)K_{g,T}\;,\quad \mbox{where}\quad
K_{g,T}=\sup \{|g(x)-g(y)|/|x-y|\;:\;x,y\in T\}\;\;.
\end{equation*}
In addition, we observe that Property $(D,d)$ implies the following simple
upper bound $K_{d,D\cap \Delta _{i}}\leq (1/u^{2})\delta $ on the Lipschitz
coefficient $K_{d,D\cap \Delta _{i}}$ of the function $d$ on the set $D\cap
\Delta _{i}$. These facts combined give 
\begin{equation*}
\mu ^{\ast }[d(D\cap \Delta _{i})]\leq \mathrm{diam}[d(D\cap \Delta
_{i})]\leq \mathrm{diam}(D\cap \Delta _{i})K_{d,D\cap \Delta _{i}}\leq
\delta ^{2}/u^{2}
\end{equation*}
whence by (\ref{bound1}) 
\begin{equation*}
\mu ^{\ast }(G)\leq 2N\max_{i=-N,N-1}\mu ^{\ast }[d(D\cap \Delta _{i})]\leq
2N\delta ^{2}/u^{2}=(2w^{2}/u^{2})/N\rightarrow 0
\end{equation*}
as $N\rightarrow \infty $ and so $\mu ^{\ast }(G)=0.\diamond $

\bigskip

The following lemma is essential in our approach to tackle the
characterization problem. We formulate it as it appears in Gikhman and
Skorohod, (1975) by combining Theorem 1 on page 263 and Theorem 4 on page
270 (see also Jacod (1985)).

\begin{lemma}
\label{charlm} Let $X(t)$ be a stochastic process with independent
increments and with $X(0)=0$. Then, for any positive number $a$, $X(t)$
admits the representation: 
\begin{eqnarray*}
X(t) &=&B(t)+\Big[c(t)+\sum_{t_{k}\leq t}\xi
_{k}+\int_{x>a}xv(t,dx)+\int_{x<-a}xv(t,dx)+\int_{0<|x|\leq a}x[v(t,dx)-\Pi
(t,dx)]\Big] \\
&=&B(t)+[c(t)\;+\;\eta (t)\;+\;T_{1,a}(t)+T_{2,a}(t)+\;U_{a}(t)] \\
&=&B(t)+c(t)+Y(t)
\end{eqnarray*}
where $\eta ,B,U_{a},T_{1,a},T_{2,a}$ are independent processes with
independent increments. The process $B$ is the zero mean continuous
component of $X$ (non--homogeneous Gaussian process) with continuous
non-decreasing variance $\sigma ^{2}(t)=\mathrm{Var}(B(t))$. The process $%
\eta $ is the deterministic time jump process, i.e. is the sum of all jumps $%
\xi _{k},$ occurred at deterministic times $t_{k}\leq t$ where the set $%
\{t_{k}\}$ is at most countable. The process $v(\Delta ,A)$ counts the
number of jumps of $X$ in a set $A$ in the interval of time $\Delta $ and $%
v(t,A)=v([0,t],A)$, where $v$ is stochastically continuous (that is $%
v(\{t\},A)=0)$. Given $A\subset R-[-a,a]$ for some $a>0$, $v(t,A)$ is a
non-homogeneous Poisson process. The measure $\Pi ((a,b),A)=\E[v((a,b),A)]$
(so $\Pi (t,A)=\E[v(t,A)]$) is its compensator. Moreover, 
\begin{equation*}
G(t,a)=\int_{0<|x|\leq a}x^{2}\Pi (t,dx)<\infty \quad \mbox{and}\quad
G(t,a)\rightarrow 0\quad \quad \mbox{as }\quad a\rightarrow 0 \;\;.
\end{equation*}
\end{lemma}

For future analysis of the processes that appear in the above representation
it is convenient to introduce the following two notations:\newline
Consider a stochastic process $\{Z=Z(s);s\in \lbrack 0,1]\}$. We say that
the process $Z$ is $f$-negligible if 
\begin{equation}
\mu \left\{ s\in \lbrack 0,1):\lim \sup_{h\downarrow 0}\;\E %
f(h^{-1/2}[Z(s+h)-Z(s)])>0\right\} \;=\;0\text{ \ .}  \label{fnegl}
\end{equation}
Next, we consider a family of stochastic processes $\{Z_{a}=Z_{a}(s);s\in
\lbrack 0,1]\}$ parameterized by $a\geq 0.$ We say that the family $%
\{Z_{a}\} $ is approximately $f$--negligible if for any real $r>0$, 
\begin{equation}
\limsup_{a\rightarrow 0}\mu \left\{ s\in \lbrack 0,1):\lim \sup_{h\downarrow
0}\;\E f(h^{-1/2}[Z_{a}(s+h)-Z_{a}(s)])>r\right\} \;=\;0\;\;.
\label{fneglapr}
\end{equation}

\bigskip

Next lemma provides some general properties about $f$--negligible processes.

\begin{lemma}
\label{fnggen} (a) If two processes $Z_{1}$ and $Z_{2}$ satisfy (\ref{fnegl}%
), and $a_{1}$ and $a_{1}$ are two real numbers, then the process $%
a_{1}Z_{1}+a_{2}Z_{2}$ is also $f$-negligible.\newline
(b) Assume that for any $a\geq 0$ the stochastic process $\{Z=Z(s);s\in
\lbrack 0,1]\}$ admits the decomposition $Z=Z_{a}+S_{a}$. If for any $a$,
the process $S_{a}$ is $f$--negligible and the family $\{Z_{a}\}$ is
approximately $f$-negligible, then $Z$ is $f$--negligible.\newline
(c) Suppose that the stochastic process $\{Z(s);s\in \lbrack 0,1]\}$
satisfies the inequality $\E f(Z(s+h)-Z(s))\leq q(s+h)-q(s)$ where $q(s)$ is
a bounded non-decreasing function. Then, the process $W$ is $f$--negligible.%
\newline
(d) Consider a family of stochastic processes $\{Z_{a}=Z_{a}(s);s\in \lbrack
0,1]\}$ parameterized by $a\geq 0,$ and suppose that $\E%
|Z_{a}(s+h)-Z_{a}(s)|^{2}\leq q_{a}(s+h)-q_{a}(s)$ where each function $%
q_{a}(s)$ is bounded, non-decreasing and $q_{a}(1)\rightarrow 0$ as $%
a\rightarrow 0$. Then, the family $\{Z_{a}\}$ is approximately $f$%
--negligible.
\end{lemma}

\textit{Proof.\/} The first and second properties are immediate consequences
of the fact that $f(x+y)\leq c_{f}[f(x)+f(y)]$ (stated in Lemma \ref{fprop})
and the additivity of the Lebesgue measure.

To prove the third property, we notice that, by Condition (\ref{fcond}) and
the condition imposed in this lemma, 
\begin{equation*}
\E f(h^{-1/2}[Z(s+h)-Z(s)])\leq Ch^{-p/2}\E f(Z(s+h)-Z(s))\leq
Ch^{-p/2}(q(s+h)-q(s))
\end{equation*}
and, since $1\leq p<2,$ it remains to apply Property \ref{pr1}.

Finally to prove Statement (d), we let $r>0$ and apply first Lemma \ref
{fprop} and then the Cauchy--Schwartz inequality to derive 
\begin{eqnarray*}
&&\mu \left\{ s\in \lbrack 0,1):\lim \sup_{h\downarrow 0}\E %
f(h^{-1/2}[Z_{a}(s+h)-Z_{a}(s)])\geq r\right\} \\
&\leq &\mu \left\{ s\in \lbrack 0,1):\lim \sup_{h\downarrow 0}[h^{-1}\E%
|Z_{a}(s+h)-Z_{a}(s)|^{2}+(h^{-1}\E|Z_{a}(s+h)-Z_{a}(s)|^{2})^{1/2}]\geq
(r/C)\right\} \;\;.
\end{eqnarray*}
Then, we apply Property \ref{pr1} along with the conditions imposed in the
part (d) of this lemma in order to bound the right hand side of the above
inequality by 
\begin{equation*}
2\mu \left\{ s\in \lbrack 0,1):\lim \sup_{h\downarrow 0}h^{-1}\E%
|Z_{a}(s+h)-Z_{a}(s)|^{2}\geq A\right\} \;\leq \;q_{a}(1)/A\rightarrow 0%
\mbox{ as }\;a\downarrow 0
\end{equation*}
(where $A=min((r/2C)^{2},(r/2C))$ ) and so the lemma follows $\diamond $

\bigskip

As one of the key steps in the proof of Theorem \ref{tm1}, we show that the
jump component is $f$--negligible which is formulated in the following lemma.

\begin{lemma}
\label{Poisson} Assume that $\E f(X(1))<\infty $. Then, the process $%
\{Y(s);s\in \lbrack 0,1]\}$ defined in Lemma \ref{charlm} satisfies (\ref
{fnegl}).
\end{lemma}

\vskip5pt \textbf{Proof}. By the property (a) of Lemma \ref{fnggen}, it is
enough to establish (\ref{fnegl}) separately for the deterministic time jump
process $\eta $ and the stochastically continuous jump process $%
T_{1,a}+T_{2,a}+U_{a}=J$, say.

We begin by analyzing the jump process $J$. By the properties (a) and (b) of
Lemma \ref{fnggen}, it is enough to show that the family $\{U_{a}\}$
satisfies (\ref{fneglapr}) and, for each $a>0$, the processes $T_{i,a}$
satisfy (\ref{fnegl}).

To show that the family $\{U_{a}\}$ is approximately $f$--negligible, we
notice that 
\begin{equation*}
\E[(U_{a}(s+h)-U_{a}(s))^{2}]=\int_{0<|x|\leq a}x^{2}\Pi
([s,s+h],dx)=G(s+h,a)-G(s,a)
\end{equation*}
where, by Lemma \ref{charlm}, for each $a>0$, the function $q_{a}(x)=G(x,a)$
is non-decreasing and $q_{a}(1)=G(1,a)\rightarrow 0$ as $a\rightarrow 0$.
Hence, (\ref{fneglapr}) is an immediate consequence of property (d) of Lemma 
\ref{fnggen}.

To finish the analysis of the stochastically continuous jump component $J$,
it is enough to show that for any $a>0$ and $i=1,2$, the process $T_{i,a}$
is $f$--negligible. Clearly, it is enough to treat only the stochastic
process 
\begin{equation*}
T_{1,a}(t)=\int_{a}^{\infty }xv(t,dx) \;\;.
\end{equation*}
By Property (c) of Lemma \ref{Glemma}, $\E f(T_{1,a}(1))<\infty $ and by
Lemma \ref{fprop}, 
\begin{equation*}
\E T_{1,a}(1)=\int_{a}^{\infty }x\Pi (t,dx)<\infty \quad \mbox{and hence}%
\quad \int_{a}^{\infty }\Pi (t,dx)<\infty \;.
\end{equation*}
Using now the week convergence approximation of the Poisson process by the
Bernoulli processes along with the Klass-Nowicki moment inequality from
Lemma \ref{KNlm}, we derive 
\begin{equation*}
c_{1}\int_{a}^{\infty }f(x)\Pi (t,dx)\;\leq \;\E f\left( \int_{a}^{\infty
}x\nu (t,dx)\right) =\E f(T_{1,a}(1))\leq c_{2}\int_{a}^{\infty }f(x)\Pi
(t,dx)\;.
\end{equation*}
Since $\E f(T_{1,a}(1))<\infty $, we note that 
\begin{equation*}
q_{a}(t)=\int_{a}^{\infty }f(x)\Pi (t,dx)\;<\;\infty
\end{equation*}
and notice that for $0\leq s<s+t\leq 1$, by Condition (\ref{fcond}) 
\begin{equation*}
\E f\left( \int_{a}^{\infty }x\nu ((s,t+s),dx)\right) \leq
c_{2}\int_{a}^{\infty }f(x)\Pi ((s,s+h),dx)\;=\;c_{2}[q_{a}(s+h)-q_{a}(s)]
\end{equation*}
and so, the process $T_{1,a}$ is $f$--negligible by Property (c) of Lemma 
\ref{fnggen}.

Now, we take care of the deterministic time jump process and notice first
that, by Property (c) of Lemma \ref{Glemma}, for any subset $A$ of the set
of points of discontinuity $\{t_{k}\}$ we have $\E f\left( \sum_{k\in A}\xi
_{k}\right) <\infty $. Furthermore, without loss of generality, we may
assume that $\E(\xi_{k})=0$ and then, by Property (d) of Lemma \ref{Glemma},
that $\xi _{k}$ are symmetric. By the Kolmogorov three series theorem and
symmetry 
\begin{equation*}
\sum_{k}P(|\xi _{k}|>1)<\infty \quad \mbox{and}\quad \sum_{k}\E(|\xi
_{k}|^{2}I_{(|\xi _{k}|\leq 1)})<\infty \;\;.
\end{equation*}
Now, for any positive $a>0$, let $Q_{a}\subseteq \{t_{k}\}$ be a finite
subset of the set of points of discontinuity such that 
\begin{equation}
\sum_{k:t_{k}\not\in Q_{a}}P(|\xi _{k}|>1)\leq a<\infty \quad \mbox{and}%
\quad \sum_{k:t_{k}\not\in Q_{a}}\E(|\xi _{k}|^{2}I_{(|\xi _{k}|\leq
1)})\leq a\;.  \label{a}
\end{equation}
We decompose the process $\eta $ into the form 
\begin{eqnarray*}
\eta (t) &=&\sum_{k:t_{k}<t,t_{k}\in Q_{a}}\xi _{k}+\sum_{k:t_{k}\leq
t,t_{k}\not\in Q_{a}}\xi _{k}I_{(|\xi _{k}|>1)}+\sum_{k:t_{k}\leq
t,t_{k}\not\in Q_{a}}\xi _{k}I_{(|\xi _{k}|\leq 1)} \\
&=&I_{1,a}+I_{2,a}+I_{3,a} \;\;.
\end{eqnarray*}
The first process $I_{1,a}$ has a finite number of jumps and obviously is $f$%
--negligible.

To analyze the second process, let $A$ be as before, a subset of the points
of discontinuity, and notice that we also have $\E f(\sum_{k\in A}\xi
_{k}I_{(|\xi _{k}|>1)})<\infty .$ Next we apply the Burkholder inequality
((1973), Theorem 15.1) and we find two constants $c_{3}$ and $c_{4}$, such
that for any subset $A\subset \{t_{k}\}$ we have 
\begin{equation}
c_{3}\E f\left( \sum_{k\in A}\xi _{k}^{2}I_{(|\xi _{k}|>1)}\right)
^{1/2}\leq \E f\left( \sum_{k\in A}\xi _{k}I_{(|\xi _{k}|>1)}\right) \leq
c_{4}\E f\left( \sum_{k\in A}\xi _{k}^{2}I_{(|\xi _{k}|>1)}\right)
^{1/2}\;\;.  \label{Burkineq}
\end{equation}

To estimate the quadratic term we apply the Klass--Nowicki inequality in
Lemma \ref{KNlm} with $X_{k}=$ $\xi _{k}^{2}$ , $B_{k}=(|\xi _{k}|>1),$ and
the function $H(x)=f(\sqrt{x}),$ $x>0$ (which obviously satisfies the
conditions of Lemma \ref{KNlm} with the power $p/2)\ $and derive 
\begin{equation}
\ c_{1}\sum_{k\in A}\E f(\xi _{k}I_{(|\xi _{k}|>1)})\leq \E f\left(
\sum_{k\in A}\xi _{k}^{2}I_{(|\xi _{k}|>1)}\right) ^{1/2}\leq
c_{2}\sum_{k\in A}\E f(\xi _{k}I_{(|\xi _{k}|>1)}) \;\;.  \label{kn}
\end{equation}
As a consequence, by (\ref{Burkineq}) and (\ref{kn}) 
\begin{equation*}
Q(t)=\sum_{k:t_{k}\leq t}\E f(\xi _{k}I_{(|\xi _{k}|>1)})<\infty \quad %
\mbox{and}\quad \E f\left( \sum_{k:s<t_{k}\leq s+h,t_{k}\not\in Q}\xi
_{k}I_{(|\xi _{k}|>1)}\right) \;\leq \;c(Q(s+h)-Q(s)) \;\;.
\end{equation*}
Therefore, the process $\{I_{2,a}(t)\}$ satisfies the conditions of Property
(c) in Lemma \ref{fnggen} and thus is $f$--negligible.

Finally, in order to treat the process $\{I_{3,a}(t)\}$ of bounded jumps, we
define the finite non-decreasing function 
\begin{equation*}
G(t)=\E\left( \sum_{k:t_{k}\leq t,t_{k}\not\in Q}\xi _{k}I_{(|\xi _{k}|\leq
1)}\right) ^{2}
\end{equation*}
and notice that 
\begin{equation*}
\E\left( \sum_{k:s<t_{k}\leq s+h,t_{k}\not\in Q}\xi _{k}I_{(|\xi _{k}|\leq
1)}\right) ^{2}=G(s+h)-G(s) \;\;.
\end{equation*}
Since by Relation (\ref{a}), $G(1)\leq a\rightarrow 0$ as $a\rightarrow 0$
it follows that the process $I_{3,a}$ satisfies the conditions of Property
(d) in Lemma (\ref{fnggen}) and therefore the family $\{I_{3,a}\}_{a\geq 0}$
is asymptotically $f$--negligible. Thus, by Property (b) in Lemma (\ref
{fnggen}) the stochastic process $\eta $ is $f$--negligible, which completes
the proof of the lemma.$\diamond $

\bigskip

The next lemma treats the Gaussian case of Theorem \ref{tm1}.

\begin{lemma}
\label{technical1} Suppose that $\{V(t);$ $t\in \lbrack 0,1]\}$ is a
stochastically continuous Gaussian process, with independent increments, $%
V(0)=0$, and there exists $\sigma \geq 0$ such that 
\begin{equation}
\lim \inf_{h\downarrow 0}\E f([V(t+h)-V(t)]/\sqrt{h})\;\geq \E f(\sigma
W)\qquad \mu -\text{almost surely} \;\;.  \label{absmean}
\end{equation}
Then, $\mathrm{Var}(V(t))\geq \sigma ^{2}t$ for all $t\in \lbrack 0,1]$ and $%
\mathrm{Var}(V(1))=\sigma ^{2}$ if and only if $\mathrm{Var}(V(t))=\sigma
^{2}t$ for all $t\in \lbrack 0,1]$.
\end{lemma}

\textit{Proof.\/} Denote by $\sigma ^{2}(t)=\mathrm{Var}(V(t))$, which is a
non-negative, continuous, non-decreasing function. First, we notice that if $%
\sigma =0,$ then the lemma is immediate.

Since $\sigma ^{2}(t)$ is non-decreasing, its derivative $(\sigma
^{2}(t))^{\prime }$ exists almost surely with respect to the Lebesgue
measure $\mu $ and to prove the lemma, it is enough to show that $\mu $
almost surely for $t\in \lbrack 0,1]$, 
\begin{equation}
\sigma ^{2}\leq (\sigma ^{2}(t))^{^{\prime }}  \label{down}
\end{equation}
Denote by $c(t)=\E V(t)$. Fix $t\in (0,1)$ such that the derivative $(\sigma
^{2}(t))^{\prime }$ exists. By Lemma \ref{funct}, there exists a positive
sequence $h_{\ast }\downarrow 0$ such that $\ h_{\ast }^{-1/2}|c(t+h_{\ast
})-c(t)|\rightarrow 0$. Then, since $f(x)$ is continuous and $|f(x)|\leq
C(|x|+x^{2})$, by the Lebesgue dominated convergence theorem, we obtain: 
\begin{eqnarray*}
&&\lim \inf_{h_{\ast }\downarrow 0}\E f(h_{\ast }^{-1/2}[V(t+h)-V(t)]) \\
&&\qquad =\lim \inf_{h_{\ast }\downarrow 0}\E f(h_{\ast
}^{-1/2}[(V(t+h)-V(t)-c(t+h_{\ast })-c(t)])=\E f(\sqrt{(\sigma
^{2}(t))^{\prime }}W)\;\;.
\end{eqnarray*}
Thus, by the lower bound in Condition (\ref{absmean}) 
\begin{equation*}
\E f(\sigma |W|)=\E f(\sigma W)\leq \E f(\sqrt{(\sigma ^{2}(t))^{\prime }}W)
\end{equation*}
for almost all $t$ which proves (\ref{down}) by Lemma \ref{Glemma}, Property
(a).

To prove the second part of the lemma we just have to notice that 
\begin{equation*}
0=\sigma ^{2}(1)-\sigma ^{2}=\int_{0}^{1}[d(\sigma ^{2}(t))-\sigma ^{2}dt]
\end{equation*}
whence, by (\ref{down}), $(\sigma ^{2}(t))^{^{\prime }}=\sigma ^{2}$, $\mu $%
-almost surely for $t\in \lbrack 0,1]$, implying that $\sigma
^{2}(t)=t\sigma ^{2}$ for all $t\in \lbrack 0,1]$.$\diamond $

\vskip5pt \textbf{Proof of Theorem \ref{tm1}.} We start from the
representation of Lemma \ref{charlm} applied to the process $\{X(t);$ $t\in
\lbrack 0,1]\}$, hence $X(t)=B(t)+c(t)+Y(t)$ for all $t\in \lbrack 0,1]$.
Since $\E|X(1)|<\infty $, by Lemma \ref{Poisson} the discontinues component,
the jump process $Y,$ satisfies $\lim \sup_{h\downarrow 0}\;\E|Y(h+s)-Y(s)|/%
\sqrt{h}=0$ almost surely with respect to the Lebesgue measure. Whence, by
condition (\ref{keycond}), the Gaussian component $\{B(t)+c(t);$ $t\in
\lbrack 0,1]\}$ satisfies (\ref{absmean}) with $\sigma =\Psi ^{-1}\E f(X(1))$%
. Denote by $\sigma ^{2}(1)=\mathrm{Var}(B(1))$. From Lemma \ref{charlm} and
Lemma \ref{Glemma} we derive 
\begin{equation}
\E f(X(1))=\E f(B(1)+c(1)+Y(1))=\E f(\sigma (1)W+c(1)+Y(1))\geq \E f(\sigma
(1)W) \;\;.  \label{ef}
\end{equation}
Moreover, by Lemma \ref{technical1}, we obtain $\sigma (1)\geq \sigma $, and
so, by the definition of $\sigma ,$ $\E f(\sigma (1)W)\geq \E f(X(1)).$ This
fact together with Relation (\ref{ef}) imply that $\mathrm{Var}(B(1))=\sigma
^{2}.$ Moreover, by the second part of Lemma \ref{technical1}, we obtain
that $\sigma ^{2}(t)=\sigma ^{2}t$ for all $t\in \lbrack 0,1]$ and, by Lemma 
\ref{Glemma}, $P(c(1)+Y(1)=0)=1,$ implying that $Y(1)$ is degenerate. Since
the process $Y(t)$ has independent increments if follows that all increments
are degenerate, which establishes (\ref{repres}). Moreover, $\E X(1)=0$
because $c(1)+Y(1)=0$ almost surely.$\diamond $

\vskip5pt \textbf{Remark and proof of Corollary \ref{tm1cor2}.\/} As it
follows from the proof of Theorem \ref{tm1}, Condition (\ref{keycond}) can
be slightly weakened to consider subsequences $h_{\ast }\rightarrow 0$ such
that the centering function $c(t)$ satisfies $h_{\ast }^{-1/2}(c(t+h_{\ast
})-c(t))\rightarrow 0$. In particular, for homogeneous stochastically
continuous processes, the centering sequence is defined by the continuous
solution of the Cauchy equation 
\begin{equation}
c(x+y)=c(x)+c(y)  \label{cauchy}
\end{equation}
implying that $c(t)=qt$, $t\in \lbrack 0,1]$. Thus, the representation (\ref
{repres}) in Corollary \ref{tm1cor2} is then immediate. Finally, $\E %
X(1)=q=0 $, which completes the proof of the corollary.$\diamond $

\vskip5pt\textbf{Proof of Proposition \ref{example}\/}. First, we choose a
positive sequence $t_{n}\downarrow 0$ such that the set $T=\{t_{n};$ $n\geq
0\}$ is independent with respect to the rational field. Then, by using Zorn
lemma, we construct the Hamel basis $B\subset R$ such that $T\subset B$. In
order to construct the function $k$ that satisfies the Cauchy equation (\ref
{cauchy}), we define it first on the set $B$ by 
\begin{equation*}
k(b)=0,\;\mbox{ if }b\notin T\;\mbox{ and }\;k(t_{i})=f^{-1}(1)\sqrt{t_{i}}%
\;,\;i=1,\ldots
\end{equation*}
Then, the solution to (\ref{cauchy}) is given by setting $c(\Sigma
r_{i}b_{i})=\Sigma r_{i}c(b_{i})$, (see for example Hardy, Littlewood and
Polya (1952)). Now, let $\{Y(t);t\geq 0\}$ be a homogeneous Poisson process
with rate $1$ and $b>0$ be such that $\E f(bY(1))=1$. Define 
\begin{equation*}
X(t)=bY(t)+k(t)-tk(1),\;t\geq 0\;\;.
\end{equation*}
Then, $\{X(t);$ $t\geq 0\}$ is a homogeneous stochastic process with
independent increments, with $X(0)=0$. Notice, that $X(1)=bY(1)$ and $\E %
f(bY(t_{n}))/\sqrt{t_{n}}\rightarrow 0$ by Lemma \ref{Poisson}, whence, by
construction, we derive 
\begin{equation*}
\E f(X(1))=\E f(bY(1))=1=\lim_{t_{n}\rightarrow 0}\E %
f(t_{n}^{-1/2}|X(t_{n})|)=\lim_{t_{n}\rightarrow 0}f(k(t_{n})/\sqrt{t_{n}}%
)=f(f^{-1}(1))=1
\end{equation*}
completing the proof of this proposition.$\diamond $

\section{Application to the Central limit theorem}

\quad This section was motivated by Theorem 3 in Dehling, Denker and Philipp
(1986). We give several applications of the characterization results from
Section 1 to extend and develop their result in several directions.

\quad The $L_{p}$ characterization of the Gaussian processes obtained in
this paper allows to avoid the traditional techniques based on the
characteristic functions in order to prove the CLT. Moreover, besides a
certain dependence condition, the additional conditions are imposed to the
moments of order $p\in[1,2)$ only. Corollary \ref{tm1cor2} is applied to
derive the following central limit theorem. Let $W$ have a standard normal
distribution and let $\Vert x\Vert _{p}=(\E|X|^{p})^{1/p}$.

\begin{theorem}
\label{tm3} Suppose that $\{X_{k};$ $k=1,2,\ldots \}$ is a strictly
stationary sequence and $p$ a fixed real, $p\in \lbrack 1,2)$. Assume $\E %
|X_{0}|^{p}<\infty $ and let $S_{n}=X_{1}+\ldots +X_{n}$, $n=1,2,\ldots $, $%
S_{0}=0$. Define the normalizing sequence $\rho _{n}=\Vert S_{n}||_{p}/\Vert
W\Vert _{p}$, $\ $ and assume that

(i) For any positive integer $k$ and real number $x$, 
\begin{equation}
\lim_{n\rightarrow \infty }\left| \E\exp \left( ixS_{n}/\rho _{n}\right)
\;-\;\left( \E\exp \left( ixS_{[n/k]}/\rho _{n}\right) \right) ^{k}\right| =0
\label{WII}
\end{equation}
(ii) $\rho _{n}\rightarrow \infty $ and there exists a positive integer $K>1$
such that $\rho _{Kn}/\rho _{n}\rightarrow \sqrt{K}$ as $n\rightarrow \infty
.$\newline
(iii) $\{(|S_{n}|/\rho _{n})^{p};n=1,2,\ldots \}$ is an uniformly integrable
family.

Then, $S_{n}/\rho _{n}\rightarrow ^{D}N(0,1)$.
\end{theorem}

\begin{corollary}
Let $\{X_{n};$ $n\geq 0\}$ be a strictly stationary sequence of integrable
random variables as in Theorem \ref{tm3} satisfying the condition (\ref{WII}%
). Let $p$ be a fixed real number $p\in \lbrack 1,2)$ and assume there is a
sequence of constants $b_{n}=\sqrt{nh(n)},$ where $h(n)$ is a function
slowly varying at $\infty $ such that the family $\{(|S_{n}|/b_{n})^{p},n%
\geq 1\}$ is uniformly integrable. Then, $\lim \Vert S_{n}\Vert _{p}/b_{n}=c$
if and only if $S_{n}/b_{n}$ converges in distribution to $N(0,\Vert W\Vert
_{p}^{2}\cdot c^{2})$.
\end{corollary}

If the second moments are finite then we immediately derive from the above
corollary:

\begin{corollary}
\label{tm3cor2}Let $\{X_{n};$ $n\geq 0\}$ be a strictly stationary sequence
of square integrable random variables satisfying the condition (\ref{WII})
and assume that $\sigma _{n}=stdev(S_{n})=\sqrt{nh(n)},$ where $h(n)$ is a
function slowly varying at $\infty .$\ Let $p$ be a fixed real number $p\in
\lbrack 1,2).$ Then, $\lim_{n\rightarrow \infty }\Vert S_{n}\Vert
_{p}/\sigma _{n}=c$ if and only if $S_{n}/\sigma _{n}$ converges in
distribution to $N(0,\Vert W\Vert _{p}^{2}\cdot c^{2}).$
\end{corollary}

Following O'Brein (1987) we say that a strictly stationary sequence $%
\{X_{k};k=1,2,\ldots \}$ is $r$-strongly-mixing sequence, if 
\begin{equation*}
\alpha _{r}(n)=\sup \left| \frac{1}{r}\left( \sum_{k=0}^{r-1}P(A\cap
B_{k})\right) -P(A)P(B)\right| \rightarrow 0\text{ as }n\rightarrow \infty
\end{equation*}
where the supremum is taken over all positive integers $m$; $A\in \mathcal{F}%
_{0}^{m}$, $\;B\in \mathcal{F}_{m+n}^{\infty }$, and $B_{k}$ is a shift of $%
B $ for $k$ steps (if $B=\{(X_{1},X_{2},\ldots )\in E\}$ for some Borel $E$,
then $B_{k}=\{(X_{k+1},X_{k+2},\ldots )\in E\}$).

It follows from Jakubowski (1993), Proposition 5.3 that $r$-strongly mixing
sequences satisfy the weak asymptotically independence condition (\ref{WII}%
). O'Brein (1987) pointed out that instantaneous functions of a stationary
Harris chain with period $d>1$ are $d$-strongly mixing and thus, by
Jakubowski (1993), they satisfy (\ref{WII}). However, they are not mixing in
a classical ergodic sense. Also, strongly mixing condition implies $r$%
--strong mixing. In particular, Theorem 3 in Dehling, Denker and Philipp
(1986) follows from Corollary \ref{tm3cor2} applied with $p=1.$\newline
\quad The regularity condition (ii) in Theorem \ref{tm3} is not easy to
check. However, using arguments similar to Jakubowski (1993) it follow that
conditions (i), (iii) and the central limit theorem $S_{n}/\rho
_{n}\rightarrow ^{D}N(0,1)$ imply (ii). Moreover, one can argue as in
Dehling, Denker and Philipp (1986) that the regularity condition can be
checked empirically, using for example the bootstrap procedure. As it is
pointed out in Peligrad (1998), the limit theorems for bootstrapped
estimators of dependent sequences require less restrictive conditions than
the corresponding limit theorems for the original sequences.

\vskip5pt \textbf{Proof of Theorem \ref{tm3}.\/} First, we derive a useful
consequence of condition (ii). We notice that, for any non--negative integer 
$j$,\quad $|\Vert S_{l+j}||_{p}/\Vert S_{l}\Vert _{p}\;-\;1|\leq \Vert
S_{l+j}-S_{l}\Vert _{p}/\Vert S_{l}\Vert _{p}\leq j\Vert X_{1}\Vert
_{p}/\Vert S_{l}\Vert _{p}\rightarrow 0$. Next, let $n=K^{r}m+j$ where $j\in
\{0,1,\ldots ,K^{r}-1\}$. Then, $[nK^{-r}]=m$, and 
\begin{equation}
\rho _{\lbrack nK^{-r}]}/\rho _{n}=\Vert S_{[nK^{-r}]}\Vert _{p}\Vert
S_{n}\Vert _{p}=(\Vert S_{m}\Vert _{p}/\Vert S_{mK^{r}}\Vert _{p})(\Vert
S_{mK^{r}}\Vert _{p}/\Vert S_{[mK^{r}+j]}\Vert _{p})\rightarrow K^{-r/2}\ 
\label{rhopr}
\end{equation}
as $n\rightarrow \infty $. Now we consider the normalized triangular array $%
S_{j}^{(n)}=S_{j}/\rho _{n}$ and observe that 
\begin{equation*}
|S_{[nK^{-r}]}^{(n)}|^{p}=|S_{[nK^{-r}]}|^{p}/\rho
_{n}^{p}=|S_{[nK^{-r}]}/\rho _{\lbrack nK^{-r}]}|^{p}|\rho _{\lbrack
nK^{-r})}/\rho _{n}|^{p}
\end{equation*}
and so, the sequence $\{|S_{[n/k]}^{(n)}|^{p};n=1,2,\ldots \}$ is uniformly
integrable by Condition (iii) of Theorem \ref{tm3} and Relation (\ref{rhopr}%
).

In order to prove this theorem it is enough to show that for any subsequence 
$(n^{\prime })\subseteq (n)$ there exists another subsequence $(n^{\prime
\prime })\subseteq (n^{\prime })$ such that $S_{n^{\prime \prime
}}^{(n^{\prime \prime })}\rightarrow ^{D}N(0,1)$. By the Helly
diagonalization technique we construct a subsequence $(n^{\prime \prime
})\subseteq (n^{\prime })$ such that $S_{[n^{\prime \prime }/k]}^{(n^{\prime
\prime })}\rightarrow ^{D}X^{(k)}$ for each positive integer $k\in N^{\prime
}$. Now, $S_{n^{\prime \prime }}^{(n^{\prime \prime })}\rightarrow
^{D}X^{(1)}$, and by Condition (i), $X^{(1)}$ is infinitely divisible (a
similar result was established in Proposition 3.1. in Samur (1984)). To
prove it, fix the integer $k$. By (i), we notice that $X_{n^{\prime \prime
},1}^{(k)}+\ldots +X_{n^{\prime \prime },k}^{(k)}\rightarrow ^{D}X^{(1)},$
where $\{X_{n^{\prime \prime },i}^{(k)};i=1\ldots ,k\}$ are $k$ independent
copies of $S_{[n^{\prime \prime }/k]}^{(n^{\prime \prime })}$. On the other
hand, it is easy to see that $X_{n^{\prime \prime },1}^{(k)}+\ldots
+X_{n^{\prime \prime },k}^{(k)}\rightarrow ^{D}X_{1}^{(k)}+\ldots
+X_{k}^{(k)}$ where $X_{1}^{(k)},\ldots ,X_{k}^{(k)}$ are independent copies
of $X^{(k)}$. By the uniqueness of the limit we obtain $X_{1}^{(k)}+\ldots
+X_{k}^{(k)}=^{D}X^{(1)},$ for any $k\geq 1$. Therefore, without loss of
generality we can take $X^{(k)}=X(1/k),$ for $k\in N^{\prime }$, where $%
\{X(t);t\geq 0\}$ is a separable homogeneous stochastic process with right
continuous sample path, with independent increments and with $X(0)=0$.
Moreover $\{X(t);t\geq 0\}$ can be assumed stochastically continuous.

Notice, $\Vert S_{[nK^{-r}]}^{(n)}\Vert _{p}=\Vert W\Vert _{p}\rho _{\lbrack
nK^{-r}]}/\rho _{n}\rightarrow K^{-r/2}\Vert W\Vert _{p}$. Since the
sequence $\{|S_{[n/k]}^{(n)}|^{p};n=1,2,\ldots \}$ is uniformly integrable,
we then derive that\quad $\lim \inf_{t_{k}\rightarrow 0}\Vert X(t_{k})\Vert
_{p}/\sqrt{t_{k}}\geq \E(|X(1)|^{p}),$ where $t_{k}=1/k$, for $k\in
N^{\prime }$ and it remains to apply Corollary \ref{tm1cor2}, which
completes the proof of the theorem.$\diamond $

\bigskip

\textbf{Acknowledgment}. The authors would like to thank Ildar Ibragimov for
pointing out relevant references. We would also like to thank Wlodek Bryc
for useful discussions on the absolute moments of sums of random variables
and on some problems solved in this paper. We are indebted to the referee
whose remarks and suggestions were very inspiring and significantly improved
the presentation of this paper.

\end{document}